# Ruled surfaces according to alternative moving frame

Burak Şahiner

*Manisa Celal Bayar University, Faculty of Arts and Sciences, Department of Mathematics, 45140 Manisa, Turkey.*

*E-mail: burak.sahiner@cbu.edu.tr*

**Abstract**

In this paper, we define a new type of ruled surface called $C-$ruled surface by using the alternative frame of a base curve. Then, we study its differential geometric properties such as striction line, distribution parameter, fundamental forms, Gaussian and mean curvatures. Moreover, we find geodesic curvatures, normal curvatures, and geodesic torsions of the base curve and the striction line of the new ruled surface. Finally, we give two related examples.

## 1. Introduction

Ruled surfaces, investigated by Gaspard Monge, are generated by the motion of a straight line along a curve called base curve in space. Despite being studied for many years, ruled surfaces are still one of the most interesting topics of differential geometry. Because of the simple structure of ruled surfaces among all surfaces, they are used in a variety of application areas including computer aided geometric design, geometric modeling, civil engineering, architecture and solid modeling [1, 2, 3, 6, 10].

The idea to correlate the straight lines called rulings which generate the ruled surface and the base curve is based on tangent surfaces. Tangent surfaces are surfaces generated by tangent vectors of a curve [13]. The study of its differential geometric properties is interesting and important in differential geometry of surfaces. Afterwards, by using the other Frenet vectors of a curve, principal normal surfaces and binormal surfaces were defined and investigated [8, 16]. Then, the method was applied for Bishop frame, type-2 Bishop frame, and rotation minimizing frame of the base curve [4, 9, 11, 12]. Thus, several special ruled surfaces were defined, and their differential geometric properties were investigated. Moreover, these special ruled surfaces were studied not only in Euclidean space but also in Minkowski space [7, 17].

Recently, Uzunoğlu et al. proposed a new frame called alternative moving frame to study the differential geometric properties of a curve in space [15]. This frame consists of three orthonormal vectors: the principal normal vector, the derivative of principal normal vector, and the Darboux vector denoted by $N$, $C$, and $W$, respectively. By using the alternative moving frame of a base curve, some new special ruled surfaces can be obtained. The ruled surfaces generated by principal normal vectors $N$ and Darboux vectors $W$ have already been examined by using Frenet frame [8]. Although these surfaces can be studied by using alternative moving frame which may be an original study, they are excluded from the scope of the present paper.

In this paper, we first give some geometric interpretations of curvatures in terms of alternative moving frame in Preliminaries section. Then we define a completely new and original ruled surface called $C-$ruled surface generated by the vectors $C$ of the alternative moving frame of the base curve. We study the differential geometric properties such as striction line, distribution parameter, fundamental forms, Gaussian and mean curvatures of the $C-$ruled surface. Also, we find the geodesic curvature, normal curvature, and geodesic torsion of the base curve and the striction line of the surface. Finally, we give two examples of $C-$ruled surfaces and investigate the differential geometric properties of them.

## 2. Preliminaries

In this section, we give some basic well-known concepts of a curve and a ruled surface in Euclidean 3-space. We also give some relations between Frenet frame and alternative moving frame of a curve and expose some geometric interpretations of curvatures in terms of alternative moving frame.

Let $\alpha(s)$ be a regular curve parameterized by arc-length parameter in Euclidean 3-space, namely, $\langle \alpha'(s), \alpha'(s) \rangle = 1$. The well-known Frenet frame $\{T, N, B\}$ along the curve $\alpha$ is defined as follows

$$T = \alpha', \ N = (1/\kappa)T', \ B = T \times N,$$

where the prime indicates the derivative with respect to $s$, $T$ is unit tangent vector field, $N$ is unit principal normal vector field, $B$ is unit binormal vector field, $\kappa = \|T'\|$ is called curvature. The Frenet derivative formulae can also be given as

$$T' = \kappa N, \ N' = -\kappa T + \tau B, \ B' = -\tau N,$$

where $\tau = -\langle B', N \rangle$ is called torsion [14]. As geometric interpretations of curvature and torsion, it can be pointed out that the curvature $\kappa$ measures the amount of deviation of the curve from the tangent line and the torsion $\tau$ measures the amount of deviation of the curve from the osculating plane spanned by the vectors $T$ and $N$. If $\kappa = 0$ (resp. $\tau = 0$), then the curve $\alpha$ is a straight line (resp. a planar curve) [5, 13]. If the ratio of the non-zero torsion to the non-zero curvature is constant, that is, $\dfrac{\tau}{\kappa}(s) = c$, where $c$ is a constant, then $\alpha$ is called helix [14].

Similar to the Frenet frame, alternative moving frame $\{N, C, W\}$ can be defined by using the unit principal normal vector $N$, the unit vector in the direction of the derivative of principal normal vector, i.e., $C = \dfrac{N'}{\|N'\|}$, and the unit vector in the direction of the Darboux vector of the Frenet frame $W = \dfrac{\tau T + \kappa B}{\sqrt{\kappa^2 + \tau^2}}$ [15]. The derivative formulas of the alternative moving frame can be given as

$$\begin{bmatrix} N' \\ C' \\ W' \end{bmatrix} = \begin{bmatrix} 0 & f & 0 \\ -f & 0 & g \\ 0 & -g & 0 \end{bmatrix} \begin{bmatrix} N \\ C \\ W \end{bmatrix} \qquad (2.1)$$

where the prime indicates the derivative with respect to $s$ which is the arc-length parameter of the curve $\alpha$, $f$ and $g$ are curvatures of the curve $\alpha$ in terms of alternative moving frame [15]. $f$ and $g$ can be called first and second alternative curvatures of alternative moving frame, respectively. They can be expressed in terms of the curvatures of Frenet frame as [15]

$$f = \sqrt{\kappa^2 + \tau^2} \qquad (2.2)$$

and

$$g = \dfrac{\kappa^2}{\kappa^2 + \tau^2} \left( \dfrac{\tau}{\kappa} \right)'. \qquad (2.3)$$

By using Eqs. (2.2), (2.3) and the geometric meaning of the curvatures of Frenet frame, we can give the geometric interpretations of the alternative curvatures in the following lemmas without any proofs.

**Lemma 2.1.** Let $\alpha$ be a curve in Euclidean 3-space and $f$ be the first alternative curvature of the curve $\alpha$. The curve $\alpha$ is a straight line if and only if $f = 0$.

**Lemma 2.2.** Let $\alpha$ be a curve in Euclidean 3-space with alternative curvatures $f \neq 0$ and $g$. The curve $\alpha$ is a helix if and only if $g = 0$.

From Eqs. (2.2) and (2.3), the curvatures of Frenet frame can be expressed in terms of the alternative curvatures. These relations can be given in the following lemma.

**Lemma 2.3.** Let $\alpha$ be a curve in Euclidean 3-space, $\{\kappa, \tau\}$ and $\{f, g\}$ be curvatures of the Frenet frame and the alternative moving frame, respectively. The curvatures of Frenet frame in terms of alternative curvatures can be given as

$$\kappa(s) = f(s) \cos\left( \int g(s) ds \right) \text{ and } \tau(s) = f(s) \sin\left( \int g(s) ds \right).$$

**Proof.** Let the ratio $\dfrac{\tau(s)}{\kappa(s)} = h(s)$. The second alternative curvature can be rewritten as $g(s) = \dfrac{h'(s)}{1 + h^2(s)}$. By integrating the last equation, we have $\int g(s) ds = \arctan(h(s))$. So, we get $h(s) = \tan\left( \int g(s) ds \right)$. By putting $\dfrac{\tau(s)}{\kappa(s)}$ instead of $h(s)$, using trigonometric relations and the first alternative moving frame $f = \sqrt{\kappa^2 + \tau^2}$, we have

$$\kappa(s) = f(s)\cos\left(\int g(s)ds\right) \text{ and } \tau(s) = f(s)\sin\left(\int g(s)ds\right). \qquad \square$$

The principal normal vector $N$ is a common vector of Frenet frame and alternative moving frame. By using $W = \dfrac{\tau T + \kappa B}{\sqrt{\kappa^2 + \tau^2}}$ and $C = W \times N$, we get $C = \dfrac{-\kappa T + \tau B}{\sqrt{\kappa^2 + \tau^2}}$. Thus, we have the alternative moving frame in terms of Frenet frame and curvatures. By using these relations and Lemma 2.3, we can give the following lemma without any proof.

**Lemma 2.4.** Let $\alpha$ be a curve in Euclidean 3-space, $\{T, N, B\}$ be Frenet frame and $\{N, C, W, f, g\}$ be alternative moving frame apparatus of the curve $\alpha$. The Frenet vectors $T$ and $B$ can be given in terms of alternative moving frame apparatus as

$$T(s) = -\cos\left(\int g(s)ds\right)C(s) + \sin\left(\int g(s)ds\right)W(s)$$

and

$$B(s) = \sin\left(\int g(s)ds\right)C(s) + \cos\left(\int g(s)ds\right)W(s).$$

**Remark 2.5.** For a planar curve, i.e., $\tau = 0$, the unit Darboux vector $W$ equals to the unit binormal vector $B$. It can be seen from Lemma 2.4 that it is possible if and only if $\sin\left(\int g(s)ds\right) = 0$.

Now we give some differential geometric properties of ruled surfaces.

A ruled surface $S$ is a surface generated by a line $L$ moving in space. Let $\alpha(s)$ be a regular curve in Euclidean 3-space and $X(s)$ be a direction vector of the line $L$. The parametric representation of the ruled surface $S$ can be given as

$$\varphi(s, v) = \alpha(s) + v X(s),$$

where $\alpha(s)$ is called the base curve or the directrix and $X(s)$ which are various positions of the line $L$ are called rulings [13]. The striction line and the distribution parameter of the ruled surface $S$ can be given as

$$\alpha^*(s) = \alpha(s) - \frac{\langle T(s), X'(s)\rangle}{\|X'(s)\|^2} X(s) \qquad (2.4)$$

and

$$P_X = \frac{\det(T, X, X')}{\|X'\|^2} \qquad (2.5)$$

where $T$ is unit tangent vector field of the base curve $\alpha$ [5, 13]. $P_X = 0$ if and only if the ruled surface $S$ is developable.

The unit normal vector $n$ of a surface can be defined by

$$n = \frac{\varphi_s \times \varphi_v}{\|\varphi_s \times \varphi_v\|}, \tag{2.6}$$

where $\varphi_s$ and $\varphi_v$ are partial derivatives of $\varphi$ with respect to $s$ and $v$, respectively.

The first fundamental form of a surface can be given as [14]

$$I = E\,ds^2 + 2F\,dsdv + G\,dv^2,$$

where

$$E = \langle \varphi_s, \varphi_s \rangle,\ F = \langle \varphi_s, \varphi_v \rangle \text{ and } G = \langle \varphi_v, \varphi_v \rangle. \tag{2.7}$$

The second fundamental form of a surface can be given as [14]

$$II = L\,ds^2 + 2M\,dsdv + N\,dv^2,$$

where

$$L = \frac{\det(\varphi_s, \varphi_v, \varphi_{ss})}{\sqrt{EG - F^2}},\ M = \frac{\det(\varphi_s, \varphi_v, \varphi_{sv})}{\sqrt{EG - F^2}} \text{ and } N = \frac{\det(\varphi_s, \varphi_v, \varphi_{vv})}{\sqrt{EG - F^2}}. \tag{2.8}$$

The Gaussian and mean curvatures of a surface can be given by

$$K = \frac{LN - M^2}{EG - F^2} \text{ and } H = \frac{EN - 2FM + GL}{2(EG - F^2)}, \tag{2.9}$$

respectively. A surface is minimal if and only if its mean curvature vanishes, i.e., $H = 0$ [14].

Let $\beta$ be a regular curve lying on a surface. The geodesic curvature, the normal curvature and the geodesic torsion of the curve $\beta$ can be given as

$$\kappa_g = \langle n_\beta \times T_\beta, T_\beta' \rangle,\ \kappa_n = \langle \beta'', n_\beta \rangle,\ \tau_g = \langle n_\beta \times n_\beta', T_\beta' \rangle, \tag{2.10}$$

respectively, where $n_\beta$ is the unit normal vector of the surface along the curve $\beta$ and $T_\beta$ is the unit tangent vector of $\beta$. By the aid of the curvatures, some special curves lying on a surface can be defined as follows.

**Definition 2.6.** ([13]) Let $\beta$ be a regular curve lying on a surface.

**(i)** The curve $\beta$ is a geodesic curve if only if $\kappa_g = 0$.

**(ii)** The curve $\beta$ is an asymptotic line if and only if $\kappa_n = 0$.

**(iii)** The curve $\beta$ is a principal line if and only if $\tau_g = 0$.

## 3. The differential geometry of $C$ – ruled surfaces

In this section, we define $C$ – ruled surface and examine its differential geometric properties.

Let $\alpha$ be a curve in Euclidean space, $\{N,C,W\}$ and $\{f,g\}$ denote alternative moving frame and alternative curvatures, respectively. A $C$ – ruled surface whose base curve is $\alpha$ can be defined as follows.

**Definition 3.1.** Let $\alpha$ be a curve except a straight line in Euclidean 3-space. The $C$ – ruled surface can be given by the following parameterization as

$$\varphi(s,v) = \alpha(s) + v\, C(s). \tag{3.1}$$

where $C$ is the unit vector of alternative moving frame of $\alpha$.

**Remark 3.2.** In Definition 3.1, we give a restriction such that the base curve $\alpha$ cannot be a straight line. By considering Eq. (2.2), the restriction means that the first alternative curvature of the curve $\alpha$ cannot equal to zero, that is, $f \neq 0$. If the curve $\alpha$ is a straight line, then the $C$ – ruled surface degenerates to the base curve $\alpha$.

By using Eqs. (2.4), (2.1), and Lemma 2.4, the striction line of a $C$ – ruled surface according to the base curve can be given as

$$\alpha^* = \alpha - \frac{g \sin\left(\int g\, ds\right)}{f^2 + g^2} C. \tag{3.2}$$

From Eqs. (3.2) and (2.3), we can give the following corollary.

**Corollary 3.3.** Let $S$ be a $C$ – ruled surface given by parametric representation $\varphi(s,v) = \alpha(s) + v\, C(s)$, where $\alpha$ is a curve except a straight line. The striction line of $S$ is also the base curve if and only if the base curve is a helix or a planar curve.

**Proof.** Let the striction line of $S$ be also the base curve. From Eq. (3.2), we have either $g=0$ or $\sin\left(\int g\, ds\right) = 0$. From Lemma 2.2 and Remark 2.5, the base curve is a helix or a planar curve. The converse can be shown similarly by using Lemma 2.2 and Remark 2.5.

By using Eqs. (2.5), (2.1), and Lemma 2.4, the distribution parameter of a $C$ – ruled surface can be obtained by

$$P = \frac{f}{f^2 + g^2} \sin\left(\int g\, ds\right). \tag{3.3}$$

From Eq. (3.3), we can give the following corollary.

**Corollary 3.4.** Let $S$ be a $C$ – ruled surface given by parametric representation $\varphi(s,v) = \alpha(s) + v\, C(s)$. $S$ is developable if and only if the base curve $\alpha$ is a planar curve.

**Proof.** The proof is clear from Eq. (3.3) and Remark 2.5.

Now, we can give the first and second fundamental forms of a $C$-ruled surface.

**Theorem 3.5.** The first fundamental form of a $C$-ruled surface can be given as

$$I = Eds^2 + 2Fdsdv + Gdv^2,$$

where $E = v^2 f^2 + \cos^2\left(\int g\, ds\right) + \left(\sin\left(\int g\, ds\right) + vg\right)^2$, $F = -\cos\left(\int g\, ds\right)$, and $G = 1$.

**Proof.** By using Eq. (2.1) and Lemma 2.4, the partial derivatives of the $C$-ruled surface with respect to $s$ and $v$ can be obtained as follows

$$\varphi_s = -vfN - \cos\left(\int g\, ds\right)C + \left(\sin\left(\int g\, ds\right) + vg\right)W, \quad \varphi_v = C. \tag{3.4}$$

By using Eqs. (2.7), we have

$$E = v^2 f^2 + \cos^2\left(\int g\, ds\right) + \left(\sin\left(\int g\, ds\right) + vg\right)^2, \quad F = -\cos\left(\int g\, ds\right), \quad G = 1. \qquad \square$$

**Theorem 3.6.** The second fundamental form of a $C$-ruled surface can be given as

$$II = Lds^2 + 2M\, dsdv$$

where

$$L = \frac{v^2(f'g - f g') + \sin\left(\int g\, ds\right)\left(v f' - f \cos\left(\int g\, ds\right)\right) - vfg}{\sqrt{v^2 f^2 + \left(\sin\left(\int g\, ds\right) + vg\right)^2}}, \quad M = \frac{f \sin\left(\int g\, ds\right)}{\sqrt{v^2 f^2 + \left(\sin\left(\int g\, ds\right) + vg\right)^2}}.$$

**Proof.** By using Eq. (2.1) and Remark 2.5, the second partial derivatives of the $C$-ruled surface can be obtained as follows

$$\varphi_{ss} = \left(-v f' + f \cos\left(\int g\, ds\right)\right)N - v(f^2 + g^2)C + v g' W, \quad \varphi_{vv} = 0, \text{ and } \varphi_{sv} = -f N + g W.$$

By calculating $\det(\varphi_s, \varphi_v, \varphi_{ss})$, $\det(\varphi_s, \varphi_v, \varphi_{vv})$, $\det(\varphi_s, \varphi_v, \varphi_{sv})$ and using Eqs. (2.8), we get

$$L = \frac{v^2(f'g - f g') + \sin\left(\int g\, ds\right)\left(v f' - f \cos\left(\int g\, ds\right)\right) - vfg}{\sqrt{v^2 f^2 + \left(\sin\left(\int g\, ds\right) + vg\right)^2}},$$

$$M = \frac{f \sin\left(\int g\, ds\right)}{\sqrt{v^2 f^2 + \left(\sin\left(\int g\, ds\right) + vg\right)^2}}, \text{ and } N = 0. \qquad \square$$

Since we have the coefficients $E, F, G, L, M, N$, we can give the Gaussian and mean curvatures of a $C$-ruled surface by using Eqs. (2.9).

**Corollary 3.7.** The Gaussian and mean curvatures of a $C$–ruled surface can be given as follows:

$$K = -\frac{f^2 \sin^2\left(\int g\, ds\right)}{\left(v^2 f^2 + \left(\sin\left(\int g\, ds\right) + vg\right)^2\right)^2},$$

and

$$H = \frac{v^2(f'g - f g') + v f' \sin\left(\int g\, ds\right) + f \sin\left(\int g\, ds\right)\cos\left(\int g\, ds\right) - vfg}{2\left(v^2 f^2 + \left(\sin\left(\int g\, ds\right) + vg\right)^2\right)^{3/2}}, \qquad (3.6)$$

respectively.

By using Corollary (3.7), we can give the following corollary.

**Corollary 3.8.** Let $S$ be a $C$–ruled surface given by parametric representation $\varphi(s,v) = \alpha(s) + v\, C(s)$. If the base curve $\alpha$ is a planar curve, then $M$ is a minimal surface.

**Proof.** Let the base curve $\alpha$ of the ruled surface $S$ be a planar curve. From Eq. (2.3) and Remark 2.5, we have $g = 0$ and $\sin\left(\int g(s)ds\right) = 0$. By substituting these values into Eq. (3.6), we get $H = 0$ which means the ruled surface is minimal. □

**Theorem 3.9.** Let $S$ be a $C$–ruled surface. The unit normal vector of the surface $S$ denoted by $n$ can be obtained as

$$n = \frac{1}{\sqrt{v^2 f^2 + \left(\sin\left(\int g\, ds\right) + vg\right)^2}} \left(-\left(\sin\left(\int g\, ds\right) + vg\right)N - vf\, W\right).$$

**Proof.** The proof is clear from Eqs. (2.6) and (3.4).

**Theorem 3.10.** Let $S$ be a $C$–ruled surface and $\alpha$ be the base curve of $S$. The geodesic curvature, the normal curvature, and the geodesic torsion of the base curve $\alpha$ can be given as

$$\kappa_g = 0,\ \kappa_n = -f \cos\left(\int g\, ds\right),\ \tau_g = 0,$$

respectively.

**Proof.** The surface normal along the base curve $\alpha$ of the ruled surface $S$ can be obtained as $n_\alpha = -N$. By using Eqs. (2.10) and Lemma 2.4, the curvatures of the base curve can be calculated directly.

From Theorem 3.10, we can give the following corollary.

**Corollary 3.11.** Let $S$ be a $C$–ruled surface and $\alpha$ be base curve of $S$.

**(i)** The curve $\alpha$ is both a geodesic curve and a principal line.

**(ii)** The curve $\alpha$ cannot be an asymptotic line.

**Proof. (i)** Since the geodesic curvature and the geodesic torsion of the base curve vanish, from Definition 2.6, the base curve is both a geodesic curve and a principal line.

**(ii)** In the definition of the $C$–ruled surface, we assume that the base curve $\alpha$ cannot be a straight line, that is, $f \neq 0$. So, the normal curvature $\kappa_n$ of the base curve cannot equal to zero which means that the curve $\alpha$ cannot be an asymptotic line. $\square$

Now, we investigate the curvatures of striction line of a $C$–ruled surface.

**Theorem 3.12.** Let $S$ be a $C$–ruled surface and $\alpha^*$ be striction line of $S$. The geodesic curvature, the normal curvature, and the geodesic torsion of the striction line $\alpha^*$ can be given as

$$\kappa_g^* = \frac{f g \cos^2\left(\int g\, ds\right)}{\sqrt{f^2 + g^2}}, \quad \kappa_n^* = -\frac{f^2 \cos\left(\int g\, ds\right)}{\sqrt{f^2 + g^2}}, \quad \tau_g^* = f g \cos^2\left(\int g\, ds\right),$$

respectively.

**Proof.** The surface normal along the striction line $\alpha^*$ of the ruled surface $S$ can be obtained as $n_{\alpha^*} = -\frac{f}{\sqrt{f^2 + g^2}} N + \frac{g}{\sqrt{f^2 + g^2}} W$. By using Eqs. (2.10) and Lemma 2.4, the curvatures of the striction line can be calculated directly.

From Theorem 3.12, we can give the following corollary.

**Corollary 3.13.** Let $S$ be a $C$–ruled surface, $\alpha$ and $\alpha^*$ be the base curve and the striction line of $S$, respectively.

**(i)** The striction line $\alpha^*$ is a geodesic curve if and only if the base curve $\alpha$ is a helix.

**(ii)** The striction line $\alpha^*$ cannot be an asymptotic line.

**(iii)** The striction line $\alpha^*$ is a principal line if and only if the base curve $\alpha$ is a helix.

**Proof. (i)** Let $\alpha^*$ be a geodesic curve. From Definition 2.6, $\kappa_g^* = \frac{f g \cos^2\left(\int g\, ds\right)}{\sqrt{f^2 + g^2}} = 0$. So, we get $g = 0$. From Lemma 2.2, the base curve $\alpha$ is a helix. The converse can be shown easily.

**(ii)** Let $\alpha^*$ be an asymptotic line. From Definition 2.6, $\kappa_n^* = -\dfrac{f^2 \cos\left(\int g\, ds\right)}{\sqrt{f^2 + g^2}} = 0$. Since $f \neq 0$, we have $\kappa_n^* \neq 0$ which means that the striction line $\alpha^*$ cannot be an asymptotic line.

**(iii)** Let $\alpha^*$ be a principal line. Then, $\tau_g^* = f\, g \cos^2\left(\int g\, ds\right) = 0$. Thus, we have $g = 0$ which means that the base curve $\alpha$ is a helix. The converse can be shown directly by using Eqs. (2.2), (2.3) and (2.10).

## 4. Examples

In this section, we give two examples of $C$–ruled surfaces and obtain their differential geometric properties.

**Example 4.1.** Let a space curve $\alpha(s) = \left(\cos\dfrac{s}{\sqrt{2}}, \sin\dfrac{s}{\sqrt{2}}, \dfrac{s}{\sqrt{2}}\right)$ be given. The alternative moving frame and alternative curvatures of the curve $\alpha$ can be obtained as

$$N(s) = \left(-\cos\dfrac{s}{\sqrt{2}}, -\sin\dfrac{s}{\sqrt{2}}, 0\right),\ C(s) = \left(\sin\dfrac{s}{\sqrt{2}}, -\cos\dfrac{s}{\sqrt{2}}, 0\right),\ W(s) = (0,0,1),$$

$$f = \dfrac{1}{\sqrt{2}},\ g = 0.$$

The $C$–ruled surface whose base curve is $\alpha(s)$ and rulings are $C(s)$ can be written as

$$\varphi(s,v) = \alpha(s) + v\, C(s) = \left(\cos\dfrac{s}{\sqrt{2}} + v\sin\dfrac{s}{\sqrt{2}}, \sin\dfrac{s}{\sqrt{2}} - v\cos\dfrac{s}{\sqrt{2}}, \dfrac{s}{\sqrt{2}}\right).$$

From Eq. (3.2), the base curve of the $C$–ruled surface is also the striction line. By using Eq. (3.3), the distribution parameter of the $C$–ruled surface can be obtained as $P = 1$. Since $P \neq 0$, the $C$–ruled surface is not a developable surface. The first-order partial derivatives of $\varphi(s,v)$ can be found as

$$\varphi_s = \left(-\dfrac{1}{\sqrt{2}}\sin\dfrac{s}{\sqrt{2}} + \dfrac{v}{\sqrt{2}}\cos\dfrac{s}{\sqrt{2}}, \dfrac{1}{\sqrt{2}}\cos\dfrac{s}{\sqrt{2}} + \dfrac{v}{\sqrt{2}}\sin\dfrac{s}{\sqrt{2}}, \dfrac{1}{\sqrt{2}}\right),$$

$$\varphi_v = \left(\sin\dfrac{s}{\sqrt{2}}, -\cos\dfrac{s}{\sqrt{2}}, 0\right).$$

By using Theorem 3.5, the coefficients of first fundamental form can be obtained as

$$E = \dfrac{v^2}{2} + 1,\ F = -\dfrac{1}{\sqrt{2}},\ G = 1.$$

Thus, the first fundamental form of the $C$–ruled surface can be written as

$$I = \left(\frac{v^2}{2}+1\right)ds^2 - \sqrt{2}\,ds\,dv + dv^2.$$

The second-order partial derivatives of $\varphi(s,v)$ can be found as

$$\varphi_{ss} = \left(-\frac{1}{2}\cos\frac{s}{\sqrt{2}} - \frac{v}{2}\sin\frac{s}{\sqrt{2}}, -\frac{1}{2}\sin\frac{s}{\sqrt{2}} + \frac{v}{2}\cos\frac{s}{\sqrt{2}}, 0\right),$$

$$\varphi_{vv} = (0,0,0),$$

$$\varphi_{sv} = \left(\frac{1}{\sqrt{2}}\cos\frac{s}{\sqrt{2}}, \frac{1}{\sqrt{2}}\sin\frac{s}{\sqrt{2}}, 0\right).$$

By using Theorem 3.6, the coefficients of second fundamental form can be obtained as

$$L = \frac{-1}{2\sqrt{v^2+1}},\ M = \frac{1}{\sqrt{2v^2+2}},\ N = 0.$$

Thus, the second fundamental form of the $C-$ruled surface can be written as

$$II = \left(\frac{-1}{2\sqrt{v^2+1}}\right)ds^2 + \frac{\sqrt{2}}{\sqrt{v^2+1}}ds\,dv.$$

By using Corollary 3.7, the Gaussian curvature and the mean curvature of the $C-$ruled surface can be obtained as

$$K = -\frac{1}{(v^2+1)^2} \text{ and } H = \frac{1}{2(v^2+1)^{3/2}},$$

respectively.

From Theorem 3.9, the unit normal vector of the $C-$ruled surface can be found as

$$n = \frac{1}{\sqrt{v^2+1}}\left(\cos\frac{s}{\sqrt{2}}, \sin\frac{s}{\sqrt{2}}, -v\right).$$

For $v=0$, we have the unit normal vector of the $C-$ruled surface along the base curve $\alpha$ as

$$n_\alpha = \left(\cos\frac{s}{\sqrt{2}}, \sin\frac{s}{\sqrt{2}}, 0\right).$$

From Theorem 3.10, the geodesic curvature, the normal curvature and the geodesic torsion of the base curve of the $C-$ruled surface can be obtained as

$$\kappa_g = 0,\ \kappa_n = -\frac{1}{2},\ \tau_g = 0,$$

respectively. From Corollary 3.11, the base curve is both a geodesic curve and a principal line, but not an asymptotic line.

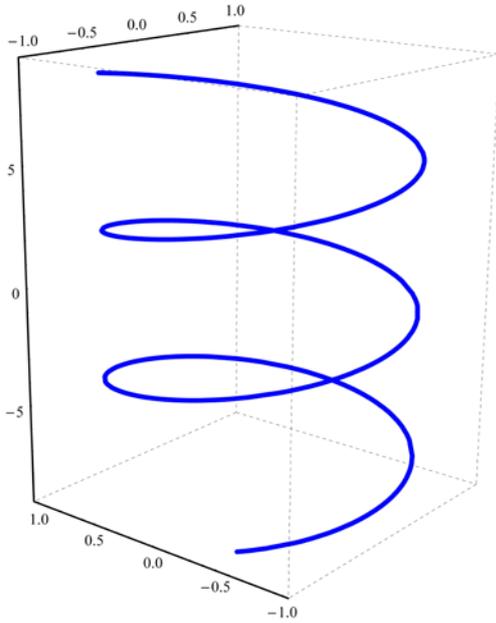
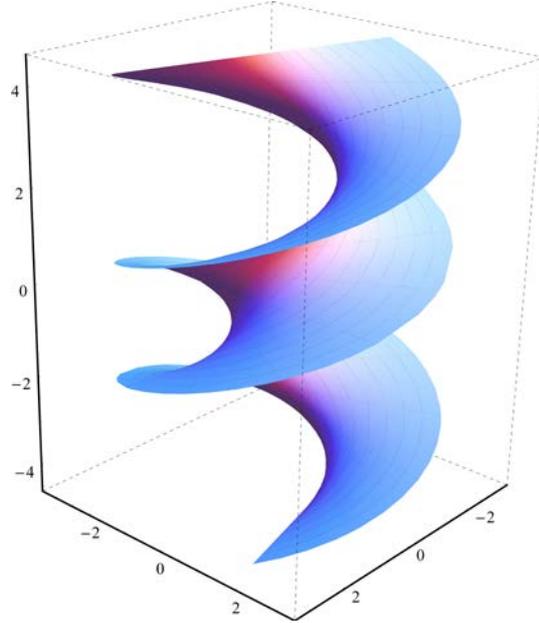

**Figure 1.a** Given curve $\alpha$.  **Figure 1.b** The $C$-ruled surface generated by the curve $\alpha$.

**Example 4.2.** Let a curve in Euclidean 3-space be given as the following parameterization

$$\alpha(s) = \left( \frac{3}{\sqrt{2}} \sin(\sqrt{2}s) \cos s - 2\sin s \cos(\sqrt{2}s), \frac{3}{\sqrt{2}} \cos(\sqrt{2}s) \cos s + 2\sin s \sin(\sqrt{2}s), -\frac{1}{\sqrt{2}} \cos s \right).$$

The alternative moving frame and alternative curvatures of the curve $\alpha$ can be obtained as

$$N(s) = \left( -\frac{1}{\sqrt{2}} \sin(\sqrt{2}\,s), -\frac{1}{\sqrt{2}} \cos(\sqrt{2}\,s), \frac{1}{\sqrt{2}} \right),$$

$$C(s) = \left( -\cos(\sqrt{2}\,s), \sin(\sqrt{2}\,s), 0 \right),$$

$$W(s) = \left( -\frac{1}{\sqrt{2}} \sin(\sqrt{2}\,s), -\frac{1}{\sqrt{2}} \cos(\sqrt{2}\,s), -\frac{1}{\sqrt{2}} \right),$$

$f = 1$, $g = -1$.

The $C$-ruled surface whose base curve is $\alpha(s)$ and rulings are $C(s)$ can be written as

$$\varphi(s,v) = \alpha(s) + v\,C(s) = \left(\frac{3}{\sqrt{2}}\sin\left(\sqrt{2}s\right)\cos s - 2\sin s\cos\left(\sqrt{2}s\right) - v\cos\left(\sqrt{2}s\right),\right.$$
$$\left.\frac{3}{\sqrt{2}}\cos\left(\sqrt{2}s\right)\cos s + 2\sin s\sin\left(\sqrt{2}s\right) + v\sin\left(\sqrt{2}s\right), -\frac{1}{\sqrt{2}}\cos s\right).$$

From Eq. (3.2), the striction line of the $C$-ruled surface can be found as

$$\alpha^* = \left(-\frac{3}{2}\sin s\cos\left(\sqrt{2}s\right) + \frac{3}{\sqrt{2}}\sin\left(\sqrt{2}s\right)\cos s, \frac{3}{2}\sin s\sin\left(\sqrt{2}s\right) + \frac{3}{\sqrt{2}}\cos\left(\sqrt{2}s\right)\cos s, -\frac{1}{\sqrt{2}}\cos s\right).$$

By using Eq. (3.3), the distribution parameter of the $C$-ruled surface can be obtained as $P = -\dfrac{\sin s}{2}$. Since $P \neq 0$, the $C$-ruled surface is not developable. By finding the first-order partial derivatives of $\varphi(s,v)$ and using Theorem 3.5, the coefficients of first fundamental form can be obtained as

$$E = 2v\sin s + 2v^2 + 1,\ F = -\cos s,\ G = 1.$$

Thus, the first fundamental form of the $C$-ruled surface can be written as

$$I = \left(2v\sin s + 2v^2 + 1\right)ds^2 - 2\cos s\,ds\,dv + dv^2.$$

By finding the second-order partial derivatives of $\varphi(s,v)$ and using Theorem 3.6, the coefficients of second fundamental form can be obtained as

$$L = \frac{\cos s\,(\sin s + v)}{\sqrt{2v\sin s + 2v^2 + 1 - \cos^2 s}},\ M = -\frac{\sin s}{\sqrt{2v\sin s + 2v^2 + 1 - \cos^2 s}},\ N = 0.$$

Thus, the second fundamental form of the $C$-ruled surface can be written as

$$II = \left(\frac{\cos s\,(\sin s + v)}{\sqrt{2v\sin s + 2v^2 + 1 - \cos^2 s}}\right)ds^2 - \frac{2\sin s}{\sqrt{2v\sin s + 2v^2 + 1 - \cos^2 s}}\,ds\,dv.$$

By using Corollary 3.7, the Gaussian curvature and the mean curvature of the $C$-ruled surface can be obtained as

$$K = -\frac{1}{(v^2+1)^2} \text{ and } H = \frac{1}{2(v^2+1)^{3/2}},$$

respectively.

From Theorem 3.9, the unit normal vector of the $C$-ruled surface can be found as

$$n = \frac{\sqrt{2}}{2\sqrt{2v\sin s + 2v^2 + 1 - \cos^2 s}}\left(-\sin\left(\sqrt{2}s\right)\sin s, -\cos\left(\sqrt{2}s\right)\sin s, (\sin s + 2v)\right).$$

For $v = 0$, the unit normal vector along the base curve can be obtained as

$$n_\alpha = \frac{1}{\sqrt{2}}\left(-\sin\left(\sqrt{2}s\right), -\cos\left(\sqrt{2}s\right), 1\right).$$

From Theorem 3.10, the geodesic curvature, the normal curvature and the geodesic torsion of the base curve of the $C-$ruled surface can be obtained as

$$\kappa_g = 0, \ \kappa_n = \cos s, \ \tau_g = 0,$$

respectively. From Corollary 3.11, the base curve is both a geodesic curve and a principal line, but not an asymptotic line.

For $v = -\dfrac{\sin s}{2}$, the unit normal vector along the striction line can be obtained as

$$n_{\alpha^*} = \left(-\sin\left(\sqrt{2}s\right), -\cos\left(\sqrt{2}s\right), 0\right).$$

From Theorem 3.12, the geodesic curvature, the normal curvature and the geodesic torsion of the striction line of the $C-$ruled surface can be obtained as

$$\kappa_g = \frac{3\sqrt{2}}{7\cos^2 s + 2}, \ \kappa_n = \frac{3\cos s}{\sqrt{2}}, \ \tau_g = -\frac{18\cos s}{(7\cos^2 s + 2)^{3/2}},$$

respectively.

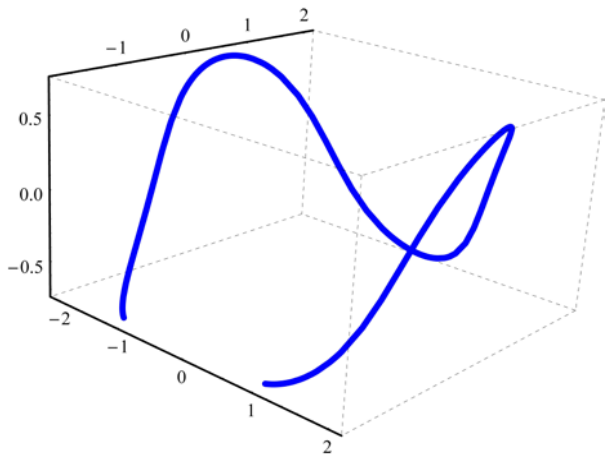 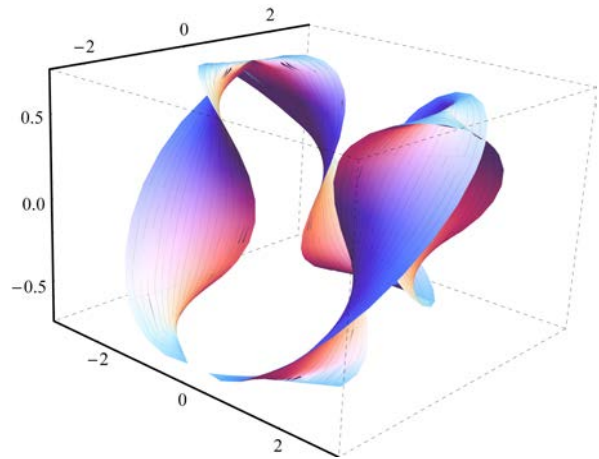

**Figure 2.a** Given curve $\alpha$.  **Figure 2.b** The $C-$ruled surface generated by the curve $\alpha$.

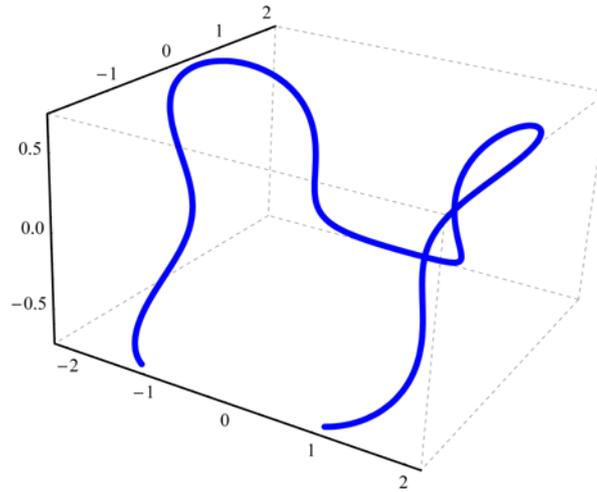

**Figure 2.c** The striction line of the $C$ – ruled surface.